\documentclass{amsart}
\usepackage{amssymb,epsfig,amsxtra, amsmath}
\usepackage[mathscr]{eucal}
\usepackage{a4}
\input{xy}
\xyoption{all}

\newtheorem{theorem}{\textbf{Theorem}}[section]
\newtheorem{proposition}[theorem]{\textbf{Proposition}}

\newtheorem{definition}[theorem]{\textbf{Definition}}
\newtheorem{lemma}[theorem]{\textbf{Lemma}}
\newtheorem{corollary}[theorem]{\textbf{Corollary}}
\newtheorem{remark}[theorem]{\textbf{Remark}}

\def\ag{\`a}

\def\p{\mathbb P^2}
\def\n{\mathbb P^{N}}

\def\s{\Sigma^n_{k,d}}

\def\g{\Gamma}

\def\bnm{\mu_{o,C}}

\def\m{\mathcal M}
\def\pa{{n -1 \choose 2}}
\def\ge{[\Gamma]}

\def\g{\Gamma}

\begin{document}
\title[On the number of moduli of plane sextics with six cusps]
{On the number of moduli of plane sextics with six cusps}

\author{ Concettina Galati }

\address{Dipartimento di Matematica, Universit\ag\, della Calabria,
Arcavacata di Rende (CS)}

\email{galati@mat.unical.it}

\thanks{}



\date{}

\dedicatory{}

\commby{}


\begin{abstract}
Let $\Sigma_{6,0}^6$ be the variety of irreducible sextics with six
cusps as singularities. Let $\Sigma\subset\Sigma_{6,0}^6$ be one of
irreducible components of $\Sigma_{6,0}^6$. Denoting by $\m_4$ the
space of moduli of smooth curves of genus $4$, we consider the
rational map $\Pi:\Sigma\dashrightarrow\m_4$ sending the general
point $\ge$ of $\Sigma$, corresponding to a plane curve
$\g\subset\p$, to the point of $\m_4$ parametrizing the
normalization curve of $\g$. The number of moduli of $\Sigma$ is, by
definition the dimension of $\Pi(\Sigma)$. We know that
$dim(\Pi(\Sigma))\leq dim(\m_4)+\rho(2,4,6)-6=7$, where
$\rho(2,4,6)$ is the Brill-Neother number of linear series of
dimension $2$ and degree $6$ on a curve of genus $4$. We prove that
both irreducible components of $\Sigma_{6,0}^6$ have number of
moduli equal to seven.
\end{abstract}

\maketitle
\section[On the number of moduli...]{Introduction}\label{introduction}
Let $\s\subset\mathbb P (H^0(\p,\mathcal{O}_{\p}(n))):=\n$, with
$N=\frac{n(n+3)}{2}$,  be the closure, in the Zariski's topology, of
the locally closed set of reduced and irreducible plane curves of
degree $n$ with $k$ cusps and $d$ nodes. Let $\Sigma\subset\s$ be an
irreducible component of the variety $\s$. Denoting by $\m_g$ the
moduli space of smooth curves of genus $g=\pa-k-d$, it is naturally
defined a rational map
$$
\Pi_{\Sigma}:\Sigma\dashrightarrow\m_g,
$$
sending the general point $[\g]\in\Sigma_0$ to the isomorphism class
of the normalization of the curve $\g\subset\p$ corresponding to
$[\g]$. We say that $\Pi_{\Sigma}$ is the \textit{moduli map of}
$\Sigma$ and we set
$$
\textit{number of moduli of $\Sigma$}:= \dim (\Pi_{\Sigma}(\Sigma)).
$$
We say that $\Sigma$ has \textit{general moduli} if $\Pi_{\Sigma}$
is dominant. Otherwise, we say that $\Sigma$ has \textit{special
moduli} or that $\Sigma$ has \textit{finite number of moduli}. By
lemma 2.2 of \cite{artmoduli}, we know that the dimension of the
general fibre of $\Pi_{\Sigma}$ is at least equal to
$$max(8,8+\rho-k),$$
where $\rho:=\rho(2,g,n)=3n-2g-6$ is the number of Brill-Noether of
linear series of degree $n$ and dimension $2$ on a smooth curve of
genus $g$. It follows that, if $\Sigma$ has the expected dimension
equal to $3n+g-1-k$ and $g\geq 2$, then
\begin{equation}
dim(\Pi_{\Sigma}(\Sigma))\leq \min(\dim(\m_g),
\dim(\m_g)+\rho-k).\label{disuguaglianza}
\end{equation}
\begin{definition}
We say that $\Sigma$ has the expected number of moduli if equality
holds in \eqref{disuguaglianza}.
\end{definition}
In particular, we expect that, if $\rho -k\leq 0$, then on the
normalization curve $C$ of the curve $\Gamma\subset\p$ corresponding
to the general point $\ge\in\Sigma$, there exists only a finite
number of linear series of degree $n$ and dimension $2$ mapping $C$
to a plane curve with nodes and $k$ cusps as singularities and
corresponding to a point of $\Sigma$, (see the proof of lemma 2.2 of
\cite{artmoduli}). For a deeper discussion and a list of known
results about the moduli problem of $\s$ we refer to sections 1 and
2 of \cite{artmoduli} and related references. In particular, in
\cite{artmoduli} we have found sufficient conditions in order that
an irreducible component $\Sigma$ of $\s$ has finite and expected
number of moduli. If $\Sigma$ verifies these conditions then
$\rho(2,n,g)\leq 0$. Finally in \cite{artmoduli} we constructed
examples of families of plane curves with nodes and cusps with
finite and expected number of moduli. In this paper we consider the
particular case of the variety $\Sigma_{6,0}^6$ of irreducible
sextics with six cusps.

It was proved by Zariski (see \cite{z2}) that $\Sigma_{6,0}^6$ has
at least two irreducible components. One of them is the parameter
space $\Sigma_1$ of the family of plane curves of equation
$$f_2^3(x_0,x_1,x_2)+f^2_3(x_0,x_1,x_2)=0,$$ where
$f_2$ and $f_3$ are homogeneous polynomials of degree two and three
respectively. The general point of $\Sigma_1$ corresponds to an
irreducible sextic with six cusps on a conic as singularities.
Moreover, $\Sigma_{6,0}^6$ contains at least one irreducible
component $\Sigma_2$ whose general element corresponds to a sextic
with six cusps not on a conic as singularities and containing in its
closure the variety $\Sigma^6_{9,0}$ of elliptic sextics with nine
cusps. Recently, A. Degtyarev has proved that $\Sigma_1$ and
$\Sigma_2$ are the unique irreducible components of
$\Sigma_{6,0}^6$, (see \cite{deg}).

The moduli number of $\Sigma_1$ and $\Sigma_2$ can not be calculated
by using the result of \cite{artmoduli}. Indeed, in this case
$\rho(2,4,6)=4>0$ and then the general element of every irreducible
component of $\Sigma^6_{6,0}$ does not verify the hypotheses of
proposition 4.1 of \cite{artmoduli}. On the contrary, it is easy to
verify that, if $\g\subset\p$ is the plane curve corresponding to
the general element of one of the irreducible components of
$\Sigma_{6,0}^6$ and $C$ is the normalization curve of $\g$, then
the map $\bnm$ is injective. But, in contrast with the nodal case,
this information is not useful in order to study the moduli problem
of $\s$, (see \cite{ser} and remark 4.2 of \cite{artmoduli}). In the
proposition \ref{subvar} and corollary \ref{subvar2}, we prove that
$\Sigma_2$ has the expected number of moduli equal to seven.
Moreover, we show that there exists a stratification
$$
\Sigma_{9,0}^6\subset \Sigma^\prime\subset\tilde\Sigma\subset\Sigma_2,\\
$$
where $\Sigma^\prime$ and $\tilde\Sigma$ are respectively
irreducible components of $\Sigma_{8,0}^6$ and $\Sigma_{7,0}^6$ with
expected number of moduli. Finally, in the corollary
\ref{seisuconica}, we prove that also $\Sigma_1$ has the expected
number of moduli by using that every element of $\Sigma_1$ is the
branch locus of a triple plane.
\section[On the number of moduli...]{On the number of moduli of
complete irreducible families of plane sextics with six
cusps}
First of all we want to find sufficient conditions in order that, if
an irreducible component $\Sigma$ of $\Sigma^{n}_{k,d}$ has the
expected number of moduli, then every irreducible component
$\Sigma^\prime$ of $\Sigma^{n}_{k^\prime,d^\prime}$, containing
$\Sigma$, has the expected number of moduli. In the corollary 4.7 of
\cite{artmoduli} we considered this problem under the hypothesis
that $\Sigma$ has the expected dimension and $\rho(2,n,g)\leq 0$.
Now we are interested to the case $\rho>0$. We need the following
local result.

Let
$$
\xymatrix{
 \mathcal D\ar[d] & = & \{(a,b,x,y)|\,y^2=x^3+ax+b\} \subset
\mathbb C^2\times \mathbb A^2\\
\mathbb C^2 & &}
$$
be the versal deformation family of an ordinary cusp (see \cite{h2}
for the definition and properties of the versal deformation family
of a plane singularity). We recall that the general curve of this
family is smooth. The locus $\Delta$ of $\mathbb C^2$ of the pairs
$(a,b)$ such that the corresponding curve is singular, has equation
$27 b^2= 4 a^3$. For $(a,b)\in\Delta$ and $(a,b)\neq (0,0)$, the
corresponding curve has a node and no other singularities, whereas
$(0,0)$ is the only point parametrizing a cuspidal curve.
\begin{lemma}[\cite{h2}, page 129.]\label{jinvariant}
Let $\mathcal G\to\mathbb
C^2$ be a two parameter family of curves of genus $g\geq
2$, whose general fibre is stable and which is locally given by
$y^2=x^3+ax+b$, with $(a,b)\in\mathbb C^2$ and let
$D\subset\mathbb C^2$ be a curve passing through $(0,0)$ and
not tangent to the axis $b=0$ at $(0,0)$. Then the $j$-invariant
of the elliptic tail of the curve which corresponds to the stable
limit of $\mathcal G_{(0,0)}$, with respect to the curve $D$, doesn't
depend on $D$. Otherwise, for every $j_0\in\mathbb C$, there
exists a curve $D_ {j_0}\subset\mathbb C^2$ passing through
$(0,0)$ and tangent to the axis $b=0$ at this point, such that the
elliptic tail of the stable reduction of $\mathcal G_{(0,0)}$ with
respect to $D_{j_0}$, has $j$-invariant equal to $j_0$.
\end{lemma}

\begin{proposition}\label{subvar}
Let $\Sigma\subset\s$, with $k<3n$, be an irreducible component of
$\s$. Let $g$ be the geometric genus of the plane curve
corresponding to the general element of $\Sigma$. Suppose that
$g\geq 2$, $\rho(2,g,n)-k\leq 0$ and $\Sigma$ has the expected
number of moduli equal to $3g-3+\rho-k$. Then, every
irreducible component $\Sigma^\prime$ of $\Sigma_{k^\prime,d^\prime}^n$,
with $k^\prime=k-1$ and $d=d^\prime$ or $k=k^\prime$ and $d^\prime=d-1$,
such that  $\Sigma\subset\Sigma^\prime$, has expected
number of moduli.
\end{proposition}
\begin{proof}
First we consider the case $k^\prime=k-1$ and $d=d^\prime$. Let
$q_1,\dots,q_k$ be the cusps of $\g$.  It is well known that, since
$k<3n$ then $\ge\in \Sigma_{k-1,d}^n$. In particular, for every
fixed cusp $q_i$ of $\g$ there exists an irreducible analytic branch
$\mathcal S_i$ of $\Sigma_{k-1,d}^n$ passing through the point $\ge$
and whose general point corresponds to a plane curve $\g^\prime$ of
degree $n$ with $d$ nodes and $k-1$ cusps specializing to the
singular points of $\g$ different from $q_i$, as $\g^\prime$
specializes to $\g$. Moreover, it is possible to prove that every
$\mathcal S_i$ is smooth at the point $\ge$, see \cite{z1} or
chapter 2 of \cite{tesi}. Let  $\Sigma^\prime$ be one of the
irreducible components of $\Sigma_{k-1,d}^n$ containing $\Sigma$.
Notice that the general element of $\Sigma^\prime$ corresponds to a
curves of genus $g^\prime=g+1$. Since
$\rho(2,g^\prime,n)-k^\prime=3n-2g-2-6-k+1=\rho(2,g,n)-k-1<0$, in
order to prove the theorem it is enough to show that the general
fibre of the moduli map
$$\Pi_{\Sigma^\prime}:\Sigma^\prime\dashrightarrow\m_{g+1}$$ has
dimension equal to eight. Let us notice that the map
$\Pi_{\Sigma^\prime}$ is not defined at the general element $\ge$ of
$\Sigma$. More precisely, let $\gamma\subset\mathcal S_i\subset
\Sigma^\prime$ be a curve passing through $\ge$ and not contained in
$\Sigma$. Let $\mathcal C\to\gamma$ be the tautological family of
plane curves parametrized by $\gamma$. Let $\mathcal C^\prime\to
\gamma$ be the family obtained from $\mathcal C\to\gamma$ by
normalizing the total space. The general fibre of $\mathcal
C^\prime\to\gamma$ is a smooth curve of genus $g+1$, while the
special fibre $\mathcal C^\prime_0:=\g^\prime$ is the partial
normalization of $\g$ obtained by smoothing all the singular points
of $\g$, except the marked cusp $q_i$. If we restrict the moduli map
$\Pi_{\Sigma^\prime}$ to $\gamma$, we get a regular map which
associates to $\ge$ the point corresponding to the stable reduction
of $\g$ with respect to the family $\mathcal C^\prime\to\gamma$,
which is the union of the normalization curve $C$ of $\g$ and an
elliptic curve, intersecting at the point $q\in C$ which maps to the
cusp $q_i\in\g$. Now, let $\mathcal G\subset\Sigma^\prime\times
\m_{g+1}$ be the graph of $\Pi_{\Sigma^\prime}$, let $\pi_1:\mathcal
G\to\Sigma^\prime$ and $\pi_2:\mathcal G\to\m_{g+1}$ be the natural
projections and let $U\subset\Sigma$ be the open set parametrizing
curves of degree $n$ and genus $g$ with exactly $k$ cusps and $d$
nodes as singularities. From what we observed before, if we denote
by $\Pi_{\Sigma^\prime}(\Sigma)$ the Zariski closure in $\m_{g+1}$
of $\pi_2\pi_1^{-1}(U)$, then $\Pi_{\Sigma^\prime}(\Sigma)$ is
contained in the divisor $\Delta_1\subset\m_{g+1}$, whose points are
isomorphism classes of reducible curves which are union of a smooth
curve of genus $g$ and an elliptic curve, meeting at a point.
Denoting by $\Pi_\Sigma:\Sigma\to\m_g$ the moduli map of $\Sigma$,
the rational map
$$\Delta_1\dashrightarrow\m_g$$ which forgets the elliptic tail,
restricts to a rational dominant map
$$q:\Pi_{\Sigma^\prime}(\Sigma)\dashrightarrow
\Pi_{\Sigma}(\Sigma).$$ The dimension of the general fibre of $q$ is
at most two. Since, by hypothesis, the dimension of the fibre of the
moduli map $\Pi_{\Sigma}$ is eight, there exists only a finite
number of $g^2_n$ on $C$, ramified at $k$ points, which maps $C$ to
a plane curve $D$ such that the associated point $[D]\in\mathbb
P^{\frac{n(n+3)}{2}}$ belongs to $\Sigma$. In particular, the set of
points $x$ of $C$ such that there is a $g^2_n$ with $k$ simple
ramification points, one of which at $x$, is finite. So, the
dimension of the general fibre of $q$ is at most one. In order to
decide if the general fibre of $q$ has dimension zero or one, we
have to understand how the $j$-invariant of the elliptic tail of the
stable reduction of $\g^\prime$ with respect the family $\mathcal
C^\prime\to\gamma$, depends on $\gamma$. If $\mathcal C\to\mathbb
C^2$ is the \'etale versal deformation family of the cusp. By
versality, for every fixed cusp $p_i$ of $\g$, there exist \'etale
neighborhoods $U\stackrel {u}\to\mathbb P^{\frac{n(n+3)}{2}}$ of
$\ge$ in $\mathbb P^{\frac{n(n+3)}{2}}$, $V\stackrel{v}\to\mathbb
C^2$ of $(0,0)$ in $\mathbb C^2$ and $U_i$ of $p_i$ in the
tautological family $\mathcal U\to\mathbb P^{\frac{n(n+3)}{2}}$ with
a morphism $f:U\to V$ such that the family $U_i\to U$ is the
pullback, with respect to $f$, of the restriction to $V$ of the
versal family. By the properties of the \'etale versal deformation
family of a plane singularity, (see \cite{dh1}), we have that
$f^{-1}((0,0))$ is an \'etale neighborhood of $\ge$ in the (smooth)
analytic branch $\Sigma_{1,0}^n$ whose general element corresponds
to an irreducible plane curves with only one cusp at a neighborhood
of the cusp $q_i$ of $\g$. So,
$\dim(f^{-1}((0,0)))=\frac{n(n+3)}{2}-2$ and the map $f$ is
surjective. Moreover, if $g$ is the restriction of $f$ at
$u^{-1}(\Sigma^\prime)$, then also $g$ is surjective. Indeed,
$$g^{-1}((0,0))=f^{-1}((0,0))\cap
u^{-1}(\Sigma^\prime)=u^{-1}((\Sigma))$$ and, since $k<3n$, then
$\dim(\Sigma)=3n+g-1-k=\dim(\Sigma^\prime)-2$ and $g$ is surjective.
By using lemma \ref{jinvariant}, it follows that the general fibre
of the natural map
$\Pi_{\Sigma^\prime}(\Sigma)\to\Pi_{\Sigma}(\Sigma)$ has dimension
exactly equal to one. Therefore, $ \dim
(\Pi_{\Sigma^\prime}(\Sigma))=\dim
(\Pi_{\Sigma}(\Sigma))+1=3g-3+\rho(2,g,n)-k+1
=3(g+1)-3+\rho(2,g+1,n)-k $ By using that
$$
\dim (\Pi_{\Sigma^\prime}(\Sigma^\prime))\geq \dim (\Pi_{\Sigma^\prime}(\Sigma))+1
=3(g+1)-3+\rho(2,g+1,n)-k+1.
$$
and by recalling that, by lemma 2.2 of \cite{artmoduli}, it is
always true that $\dim(\Pi_{\Sigma^\prime}(\Sigma^\prime))\leq
3(g+1)-3+\rho(2,g+1,n)-k+1$, the statement is proved in the case
$k^\prime=k-1$ and $d^\prime=d$.

Suppose, now, that $k=k^\prime$ and $d^\prime=d-1$. Also in this
case $\Sigma$ is not contained in the regularity domain of
$\Pi_{\Sigma^\prime}$. More precisely, if $\ge\in\Sigma$ is general,
then $\Pi_{\Sigma^\prime}(\ge)$ consists of a finite number of
points, corresponding to the isomorphism classes of the partial
normalizations of $\g$ obtained by smoothing all the singular points
of $\g$, except for a node. Then $\Pi_{\Sigma^\prime}(\Sigma)$ is
contained in the divisor $\Delta_0$ of $\m_{g+1}$ parametrizing the
isomorphism classes of the analytic curves of arithmetic genus $g+1$
with a node and no more singularities.  The natural map
$\Delta_0\dashrightarrow\m_g$ sending the general point $[C^\prime]$
of $\Delta_0$ to the isomorphism class of the normalization of
$C^\prime$, restricts to a rational dominant map
$q:\Pi_{\Sigma^\prime}(\Sigma)\dashrightarrow\Pi_{\Sigma}(\Sigma)$.
Since we suppose that $\Sigma$ has the expected number of moduli and
$\rho(2,g,n)-k\leq 0$, if $C$ is the normalization of the plane
curve corresponding to the general element of $\Sigma$, then the set
$S$ of the linear series of dimension $2$ and degree $n$ on $C$ with
$k$ simple ramification points, mapping $C$ to a plane curve $D$
such that the associated point $[D]$ in the Hilbert Scheme belongs
to $\Sigma$, is finite. We deduce that also the set $S^\prime$ of
the pairs of points $(p_1, p_2)$ of $C$, such that there is a
$g^2_n\in S$ such that the associated morphism maps $p_1$ and $p_2$
to the same point of the plane, is finite. So, also $q^{-1}([C])$ is
finite and
$\dim(\Pi_{\Sigma^\prime}(\Sigma))=\dim(\Pi_{\Sigma}(\Sigma))$. It
follows that
$$
\dim (\Pi_{\Sigma^\prime}(\Sigma^\prime))\geq 3g-3+3n-2g-6-k+1=
3(g+1)-3+3n-2(g+1)-6-k.
$$
\end{proof}
\begin{remark}
Notice that, the arguments used before to prove lemma \ref{subvar}
don't work if the dimension of the general fibre of the moduli map
of $\Sigma$ has dimension bigger than eight. Indeed, in this case,
the dimension of the general fibre of the map
$\Pi_{\Sigma^\prime}(\Sigma)\dashrightarrow\Pi_{\Sigma}(\Sigma)$ could be
bigger than one if $k^\prime=k-1$ and $d=d^\prime$, or  than zero
if $k^\prime=k$ and $d=d^\prime-1$.
\end{remark}
\begin{corollary}\label{subvar2}
There exists at least one irreducible component $\Sigma_2$ of
$\Sigma^6_{6,0}$ having the expected number of moduli equal to
$\dim(\m_4)-2$ and whose general element corresponds to a sextic
with six cusps not on a conic.
\end{corollary}
\begin{remark}
As we already observed in the previous section, $\Sigma_2$ is
the only component of $\Sigma_{6,0}^6$ parametrizing sextics with
six cusps not on a conic by \cite{deg}.
\end{remark}
\begin{proof}
Let $\Sigma_{9,0}^6$ be the variety of elliptic plane curves of
degree six with nine cusps and no more singularities. It is not
empty and irreducible, because, by the Pl\"{u}cker formulas, the
family of dual curves is $\Sigma_{0,0}^3\simeq\mathbb P^9$, which is
irreducible and not empty. Moreover, if we compose an holomorphic
map $\phi:C\to\p$ from a complex torus $C$ to a smooth plane cubic
with the natural map $\phi (C)\to{\phi (C)}^*$, where we denoted by
${\phi (C)}^*$ the dual curve of $\phi (C)$, we get a morphism from
$C$ to a plane sextic with nine cusps. Therefore, the number of
moduli of $\Sigma_{9,0}^6$ is equal of the number of moduli of
$\Sigma_{0,0}^3$, equal to one. Since $6<3n=18$, there is at least
one irreducible component $\Sigma^\prime$ of $\Sigma_{8,0}^6$
containing $\Sigma_{9,0}^6$. Let
$\Pi_{\Sigma^\prime}:\Sigma^\prime\dashrightarrow\m_2$ be the moduli
map of $\Sigma^\prime$ and let $\mathcal
G\subset\Sigma^\prime\times\m_2$ be its graph. If we denote by
$\pi_1:\mathcal G\to\Sigma^\prime$ and $\pi_2:\mathcal G\to \m_2$
the natural projection, by $U$ the open set of $\Sigma_{9,0}^6$
parametrizing cubics of genus one with nine cusps and by
$\Pi_{\Sigma^\prime}(\Sigma_{9,0}^6)$ the Zariski closure in $\m_2$
of $\pi_2\pi_1^{-1}(U)$, then, by arguing as in the first part of
the proof of the lemma \ref{subvar}, we have a dominant map
$\Pi_{\Sigma^\prime}(\Sigma_{9,0}^6)\dashrightarrow\m_1$, whose
general fibre has dimension one. We conclude that
$$
\dim(\pi_{\Sigma^{\prime}}(\Sigma^{\prime}))\geq
\dim(\pi_{\Sigma^{\prime}}(\Sigma_{9,0}^6))+1=3
$$
and so, the moduli map of $\Sigma^{\prime}$ is dominant, as one
expects, because $\rho(2,2,6)-8=18-4-6-8=0$. Let $D$ be the plane
sextic corresponding to the general point of $\Sigma^\prime$. By
Bezout theorem, the height cusps $P_1,\dots,P_8$ of $D$ don't belong
to a conic and, however we choose five cusps of $D$, no four of them
lie on a line. Then, let $C_2$ be the unique conic containing
$P_1,\dots,P_5$. There exists at least a cusp, say $P_6$, which does
not belong to $C_2$. Since $8<3n=18$, there exists a family of plane
sextics $\mathcal D\to\Delta$, whose special fibre is $D$ and whose
general fibre has a cusp at a neighborhood of every cusp of $D$
different from $P_7$ and no further singularities. By lemma
\ref{subvar}, the curve $\Delta$ is contained in an irreducible
component of $\Sigma_{7,0}^6$ with expected number of moduli. By
repeating the same argument for the general fibre of the family
$\mathcal D\to\Delta$ we get an irreducible component $\Sigma_2$ of
$\Sigma_{6,0}^6$  with the expected number of moduli and whose
general element parametrizes a sextic with six cusps not on a conic.
\end{proof}
Now we consider the irreducible component $\Sigma_1$ of
$\Sigma_{6,0}^{6}$ parametrizing plane curves of equation
$f^2_3(x_0,x_1,x_2)+f^3_2(x_0,x_1,x_2)=0$, where $f_2$ is an
homogeneous polynomial of degree two and $f_3$ is an homogeneous
polynomial of degree three. The general element of $\Sigma_1$ corresponds to an
irreducible plane curve of degree six with six cusps on a conic.
We want to show that $\Sigma_1$ has the expected number of moduli
equal to $12-3+\rho(2,4,6)-6=7=\dim(\m_4)-2$. Equivalently, we want
to show that the general fibre of the moduli map
$$\Sigma_1\dashrightarrow\m_4 $$ has dimension equal to eight.
\begin{lemma}\label{sscc}
Let $\g_2:f_2(x_0,x_1,x_2)=0$ and $\g_3:f_3(x_0,x_1,x_2)=0$ be a
smooth conic and a smooth cubic intersecting transversally. Then,
the plane curve $$\g:f^2_3(x_0,x_1,x_2)-f^3_2(x_0,x_1,x_2)=0$$ is
an irreducible sextic of genus four with six cusps at the
intersection points of $\g_2$ and $\g_3$ as singularities. The
curve $\g$ is projection of a canonical curve $C\subset\mathbb
P^3$ from a point $p\in\mathbb P^3$ which is contained in six
tangent lines to $C$. Moreover, for every point $q\in\mathbb
P^3-C$ such that the projection plane curve $\pi_q(C)$ of $C$ from
$q$ is a sextic with six cusps on a conic of equation
$g^2_3(x_0,x_1,x_2)-g^3_2(x_0,x_1,x_2)=0$, where $g_3$ and $g_2$
are two homogeneous polynomials of degree three and two
respectively, there exists a cubic surface $S_3\in |\mathcal
I_{C|\mathbb P^3}(3)|$, containing $C$, such that the plane curve
$\pi_q(C)$ is the branch locus of the projection $\pi_q:S_3\to\p$.
\end{lemma}
\begin{remark}\label{pareschi}
Notice that, by \cite{deg}, every irreducible sextic with six cusps
on a conic as singularities has equation given by
$(f_2(x_0,x_1,x_2))^3+(f_3(x_0,x_1,x_2))^2=0$, with $f_2$ and $f_3$
homogeneous polynomials of degree two and three. In order words, all
the sextics with six cusps on a conic as singularities are
parametrized by points of $\Sigma_1$. An other proof of this result
as been provided to us by G. Pareschi.
\end{remark}
\begin{proof}[Proof of lemma \ref{sscc}.]
Let $f(x_0,x_1,x_2)=f^2_3(x_0,x_1,x_2)-f^3_2(x_0,x_1,x_2)=0$ be
the equation of $\g$. From the relation $f_3(\underline x)=\pm
f_2(\underline x)\sqrt{f_2(\underline x)}$, we deduce that
$\frac{\partial f_3}{\partial x_i}(\underline x)=\pm
2\frac{\partial f_2}{\partial x_i }(\underline
x)\sqrt{f_2(\underline x)}$ and hence
\begin{equation}\label{derivativesextic}
\frac{\partial f}{\partial x_i }(\underline x)=2\frac{\partial
f_3}{\partial x_i}(\underline x)f_3(\underline x)-3{f_2(\underline
x)}^2\frac{\partial f_2 }{\partial x_i}(\underline
x)=-{f_2(\underline x)}^2\frac{\partial f_2 }{\partial
x_i}(\underline x).
\end{equation}
By using that the conic $\g_2:f_2=0$ is smooth, it follows that,
if a point $\underline x\in\g$ is singular, then $\underline
x\in\g_2$ and hence $\underline x\in\g_3\cap\g_2$. On the other
hand, always from \eqref{derivativesextic}, if $\underline
x\in\g_2\cap\g_3$, then $\underline x$ is a singular point of
$\g$. Hence, the singular locus of $\g$ coincides with
$\g_3\cap\g_2$. Let $\underline x$ be a singular point of $\g$. If
$$p_1(x,y)+\textrm{terms of degree two}=0$$ and
$$q_1(x,y)+\textrm{terms of degree}\,\,\geq\,\, \textrm{two}=0$$ are
respectively affine equations of $\g_2$ and $\g_3$ at $\underline
x$, then, the affine equation of $\g$ at $x$ is given by
$$q_1(x,y)^2-p_1(x,y)^3+\textrm{terms of degree}\,\,\geq\,\, \textrm{four}=0.$$
Since $\g_2$ and $\g_3$ intersect transversally, we have that
$q_1(x,y)$ does not divide $p_1(x,y)$ and hence $\g$ has an ordinary
cusp at $x$. Let now $\phi:C\to\g$ be the normalization of $\g$. We
recall that the cubics passing through the six cusps of $\g$ cut out
on $C$ the complete canonical series $|\omega_C|$. Since the cusps
of $\g$ is contained in the conic $\g_2\subset \p$ of equation
$f_2=0$, the lines of $\p$ cut out on $C$ a subseries $g\subset
|\omega_C|$ of dimension two of the canonical series. Moreover, if
we still denote by $C$ a canonical model of $C$ in $\mathbb P^3$,
then the linear series $g$ is cut out on $C$ in $\mathbb P^3$ from a
two dimensional family of hyperplanes passing through a point
$p\in\mathbb P^3-C$. If we project $C$ from $p$ we get a plane curve
projectively equivalent to $\g$. Since $\g$ has six cusps as
singularities, we deduce that there are six tangent lines to $C$
passing through $p$. To see that $\g$ is the branch locus of a
triple plane, let $S_3\subset\mathbb P^3$ be the cubic surface of
equation
$$
F_3(x_0,\dots,x_3)=x_3^3-3f_2(x_0,x_1,x_2)x_3+2f_3(x_0,x_1,x_2)=0.
$$
If $p=[0,0,0,1]$, then, by using Implicit Function Theorem, the
ramification locus of the projection $\pi_p:S_3\to\p$, is given by
the intersection of $S_3$ with the quadric $S_2$ of equation
$\frac{\partial F_3}{\partial x_3}=x^2_3-f_2(x_0,x_1,x_2)=0.$ Now,
if $\underline x=[x_0,x_1,x_2]\in S_3\cap S_2$, then $x_3=\pm
\sqrt{f_2(x_0,x_1,x_2)}$. By substituting in the equation of $S_3$,
we find that the branch locus of the projection $\pi_p:S_3\to\p$
coincides with the plane curve $\g$. From what we proved before, it
follows that the ramification locus of the projection map
$\pi_p:S_3\to\p$ is the normalization curve $C$ of $\g$. Finally, if
$q\in\mathbb P^3-C$ is an other point such that the plane projection
$\pi_q(C)$ is an irreducible sextic with six cusps on a conic
parametrized by a point $x_q\in\Sigma_1\subset \mathbb P^{27}$,
then, up to projective motion, we may always assume that
$q=[0:0:0:1]$ and hence, if
$g^2_3(x_0,x_1,x_2)-g^3_2(x_0,x_1,x_2)=0$ is the equation of the
plane curve $\pi_q(C)$, then $C$ is the locus of ramification of the
projection from $q$ to the plane of the cubic surface of equation
$$
x_3^3-3g_2(x_0,x_1,x_2)x_3+2g_3(x_0,x_1,x_2)=0.
$$
\end{proof}

\begin{corollary}\label{seisuconica}
The irreducible component $\Sigma_1$ of $\Sigma_{6,0}^6$
parametrizing plane curves of equation
$f^2_3(x_0,x_1,x_2)+f^3_2(x_0,x_1,x_2)=0$, where $f_2$ is an
homogeneous polynomial of degree two and $f_3$ is an homogeneous
polynomial of degree three, has the expected number of moduli
equal to $7=\dim(\m_4)+\rho(2,4,6)-6$.
\end{corollary}
\begin{proof}
Let $\ge\subset\p$ be a plane sextic of equation
$f^2_3(x_0,x_1,x_2)-f^3_2(x_0,x_1,x_2)=0$, where the conic $f_2=0$
and the cubic $f_3=0$ are smooth and they intersect transversally.
Let $C\subset\mathbb P^3$ be the normalization curve of $\g$ and
let $\mathcal S_C$ be the set of points $\underline
v=[v_0:\dots:v_3]\in\mathbb P^3$ such that there exists a cubic
surface $S_3\in |\mathcal I_{C|\mathbb P^3}(3)|$, containing $C$,
such that the curve $C$ is the ramification locus of the
projection $\pi_{\underline v}:S_3\to\p$. By the former lemma, in order to
prove that $\Sigma_1$ has the expected number of moduli, it is
enough to find a point $\ge$ of $\Sigma_1$ corresponding to an
irreducible plane sextic $\g\subset\p$ with six cusps of a conic
such that the set $\mathcal S_C$ is finite. Let $\g_2$ be the
smooth conic of equation $f_2(x_0,x_1,x_2)=x_0^2+x^2_1-x_2^2=0$
and let $\g_3$ be the smooth cubic of equation $f_3(x_0,x_1,x_2)=
x_0^3+x_0x_2^2-x_1^2x_2=0$. If $a_1,a_2$ and $a_3$ are the three
different solutions of the polynomial $x^3+x^2+x-1=0$, then $\g_2$
and $\g_3$ intersect transversally at the points
$[a_i,\sqrt{a_i},1]$, $[a_i,-\sqrt{a_i},1]$, with $i=1,2,3$. By
the former lemma, the plane sextic $\g$ of equation
$f_2^3-f_3^2=0$ is irreducible and it has six cusps at the
intersection points of $\g_2$ and $\g_3$ as singularities.
Moreover, the normalization curve $C$ of $\g$ is the canonical
curve of genus $4$ in $\mathbb P^3$ which is intersection of the
cubic surface $S_3\subset\mathbb P^3$ of equation
$$
F_3(x_0,x_1,x_2,x_3)=x_3^3+(x_0^2+x_1^2-x_2^2)x_3+x_0^3+x_0x_2^2-x_1^2x_2=0
$$ and the quadric $S_2$ of equation
$$
\frac{\partial F_3}{\partial x_3}=3x_3^2+x_0^2+x_1^2-x_2^2=0.
$$
We want to show that $\mathcal S_C$ is finite. To see this we
observe that, since $$h^0(\mathbb P^3,\mathcal I_{C|\mathbb
P^3}(2))=1\,\,\mbox{and}\,\, h^0(\mathbb P^3,\mathcal I_{C|\mathbb
P^3}(3))=5,$$ the equation of every cubic surface containing $C$ and
which is not the union of $S_2$ and an hyperplane is given by
$$
G(x_0,\dots,x_3;\beta_0,\dots,\beta_3)=F_3(x_0,x_1,x_2,x_3)+\sum_{j=0}^3\beta_j x_j
\frac{\partial F_3(x_0,x_1,x_2,x_3)}{\partial x_3 }=0,
$$
with $\beta_j\in\mathbb C$, for $i=0,\dots,3$. Now, a point
$[\underline v]=[v_0,\dots,v_3]\in \mathcal S_C$ if and only if
there exist $\beta_0,\dots,\beta_3$ such that $C$ is contained in
the intersection of $$G(x_0,\dots,x_3;\beta_0,\dots,\beta_3)=0\,\,\,\mbox{ and}\,\,\,
\frac{\partial G(x_0,\dots,x_3;\beta_0,\dots,\beta_3)}{\partial
\underline v}=0.$$ Still using that $h^0(\mathbb P^3,\mathcal
I_{C|\mathbb P^3}(2))=1$, a point $[\underline v]\in\mathbb P^3$
belongs to $\mathcal S_C$ if and only if
\begin{equation}
\frac{\partial G(x_0,\dots,x_3;\beta_0,\dots,\beta_3)}{\partial
\underline v}=\lambda\frac{\partial F_3(x_0,\dots,x_3)}{\partial
x_3}
\end{equation}
for some $\lambda\in\mathbb R-0$, or, equivalently,
\begin{equation}
\sum_{i=0}^2v_i\frac{\partial F_3}{\partial x_i}+\sum_{i=0}^3
v_i(\sum_{j=0}^3\beta_jx_j)\frac{\partial F_3}{\partial
x_3\partial
x_i}=(\lambda-\sum_{i=0}^3v_i\beta_i-v_3)\frac{\partial
F_3}{\partial x_3}.
\end{equation}
The previous equality of polynomials is equivalent to the
following bilinear system of ten equations in the
variables $v_0,\dots,v_3$ and $\beta_0,\dots,\beta_3$
\begin{equation}\label{system}
  \begin{cases}
(1+\beta_3)v_0+3\beta_0v_3=0\hfill &(x_0x_3)\\
(1+\beta_3)v_1+3\beta_1v_3=0\hfill &(x_1x_3)\\
(1+\beta_3)v_2-3\beta_2v_3=0\hfill &(x_2x_3)\\
\beta_1v_0+\beta_0v_1=0\hfill&(x_0x_1)\\
\beta_2v_0+(1-\beta_0)v_2=0\hfill&(x_0x_2)\\
(1-\beta_2)v_1+\beta_1v_2=0\hfill&(x_1x_2)\\
2\beta_1v_1-v_2=\lambda-\sum_{j=0}^3\beta_jv_j-v_3&\hfill (x_1^2)\\
-v_0+2\beta_2v_2=\lambda-\sum_{j=0}^3\beta_jv_j-v_3 &\hfill (x_2^2)\\
(3+2\beta_0)v_0=\lambda-\sum_{j=0}^3\beta_jv_j-v_3&\hfill (x_0^2)\\
2\beta_3v_3=\lambda-\sum_{j=0}^3\beta_jv_j-v_3&\hfill (x_3^2)
 \end{cases}
\end{equation}
The points of $\mathcal S_C$ are the solutions $\underline v$ of
the previous linear system, as a linear system whose coefficients
depend on $\beta_0,\dots,\beta_3$. It is easy to prove that
it has only a solution equal to $(v_0,v_1,v_2,v_3)= (0,0,0,\lambda)$ if
$\beta_0=\beta_1=\beta_2=\beta_3=0$ and it has not solutions
otherwise, (see \cite{tesi}, page 98). By the previous lemma, we conclude that the point
$[0:0:0:1]\in\mathbb P^3$ is the only point which belongs to six
tangent lines to the canonical curve $C\subset\mathbb P^3$ which
is intersection of the cubic surface of equation
$$
F_3(x_0,x_1,x_2,x_3)=x_3^3+(x_0^2+x_1^2-x_2^2)x_3+x_0^3+x_0x_2^2-x_1^2x_2=0
$$
and the quadric of equation
$$
\frac{\partial F_3}{\partial x_3}=3x_3^2+x_0^2+x_1^2-x_2^2=0.
$$
It follows that, on the normalization curve $D$ of the plane curve
$\g^\prime$ corresponding to the general point of
$\Sigma_1\subset\Sigma_{6,0}^6$ there exists only a finite number
of linear series of dimension two with six ramification points.
\end{proof}

\begin{remark}
By using the notation introduced in the proof of corollary \ref{seisuconica},
we observe that in this corollary we have proved that if $C$ is a general canonical curve
of genus four such that the set $\mathcal S_C$ is not empty, then $\mathcal S_C$
is finite. Actually, C. Ciliberto pointed out to our attention that it is possible to show,
with a very simple argument, that for every canonical curve $C$ of genus four such that $\mathcal S_C$ is not empty, we have that $\mathcal S_C$ is finite.
Finally, we observe that, by remark \ref{pareschi}, for every canonical curve $C$
of genus four, the set $\mathcal S_C$ coincides
with the set of points of $\mathbb P^3$ which are contained in six tangent lines to $C$.
\end{remark}

\subsection*{Acknowledgment}
The results of this paper are contained in my PhD-thesis. I would
like to thank my advisor C. Ciliberto for introducing me into the
subject and for providing me very useful suggestions. I have also
enjoyed and benefited from conversation with G. Pareschi and my
college M. Pacini.

{}

\end{document}